\documentclass{article}
\usepackage{amsmath} 
\usepackage{mathtools}
\usepackage{amsfonts}
\usepackage{amssymb} 
\usepackage{mathabx} 
\usepackage{graphicx}
\usepackage{fancyhdr}
\usepackage{setspace}
\usepackage{lastpage}
\usepackage[utf8]{inputenc}
\usepackage[T1]{fontenc}
\usepackage{cite}
\usepackage[title]{appendix}
\usepackage{tabulary}

\usepackage
[
letterpaper,
left=2.5cm,
right=2.5cm,
top=2.5cm,
bottom=2.5cm
]
{geometry}

\newcommand{\mytitle}{Dynamic Optimization on Quantum Hardware: Feasibility for a Process Industry Use Case}
\newcommand{\myshorttitle}{Dynamic Optimization on Quantum Hardware}
\newcommand{\myauthor}{Dennis~M.~Nenno, Adrian Caspari$^{\ast}$} 
\newcommand{\myauthorshort}{D.~M.~Nenno \& A.~Caspari}
\author{\myauthor}

\usepackage[colorlinks,linkcolor=blue,citecolor=blue,urlcolor=blue,pdftitle={\mytitle},pdfauthor={Adrian Caspari},hidelinks]{hyperref} 
\usepackage{enumerate}

\usepackage{caption} 
{
	\fancyhead[R]{\myshorttitle}
	\fancyhead[L]{\myauthorshort}
	\fancyfoot[L]{}
	\fancyfoot[R]{Page \thepage\ of \pageref*{LastPage}}
	\cfoot{}
}

\fancypagestyle{firststyle}
{ 

   \fancyhead[L]{}
	 \fancyhead[C]{}
   \fancyhead[R]{}
	 \fancyfoot[L]{\footnotesize
	 	$^{*}$Corresponding author: A.~Caspari, E-mail: \url{adrian.caspari@rwth-aachen.de} \\ Process Systems Engineering (AVT.SVT), RWTH Aachen University, 52074 Aachen, Germany   
	 }
   \fancyfoot[C]{}
	 \fancyfoot[R]{Page \thepage\ of \pageref*{LastPage}}
}

\renewenvironment{abstract}{\noindent\textbf{Abstract}\\}{}

\begin{document}

\thispagestyle{firststyle}
\begin{flushleft}\begin{Large}\textbf{\mytitle}\end{Large} \end{flushleft}
{\large \myauthor}

\vspace{0.5cm}

\begin{abstract}
	The quest for real-time dynamic optimization solutions in the process industry represents a formidable computational challenge, particularly within the realm of applications like model-predictive control, where rapid and reliable computations are critical. Conventional methods can struggle to surmount the complexities of such tasks. Quantum computing and quantum annealing emerge as \textit{avant-garde} contenders to transcend conventional computational constraints. We convert a dynamic optimization problem, {characterized by an optimization problem with a system of differential-algebraic equations embedded}, into a Quadratic Unconstrained Binary Optimization problem, enabling quantum computational approaches. The empirical findings synthesized from classical methods, simulated annealing, quantum annealing via D-Wave's quantum annealer, and hybrid solver methodologies, illuminate the intricate landscape of computational prowess essential for tackling complex and high-dimensional dynamic optimization problems. Our findings suggest that while quantum annealing is a maturing technology that currently does not outperform state-of-the-art classical solvers, continuous improvements could eventually aid in increasing efficiency within the chemical process industry.
\end{abstract}

\vspace{0.5cm}


\section{Introduction}
Dynamic optimization problems, i.e. optimization problems with differential equations embedded, appear in many application areas and are crucial for both real-time applications, such as model predictive control (MPC), and offline process optimization~{\mbox{\cite{Biegler.2010}}}. 
For instance, MPC, pivotal in process control, relies on efficiently solving dynamic optimization problems to predict and optimize the future behavior of a process ~\cite{rawlings2017model}.
However, the complexity of dynamic optimization problems often precludes their real-time solutions, in particular in those sectors, where speed is crucial, as often for MPC applications in the chemical industries~{\mbox{\cite{Caspari2019b, Caspari2019g}}}. 
In recent years, quantum computers have garnered significant attention, for their potential to revolutionize optimization problem-solving across diverse fields. While traditional approaches in dynamic optimization predominantly focus on simplifying models~{\cite{Schaefer2019a, Schulze2020}}, approximating solutions~{\cite{Kadam2004, Vaupel2020}}, and enhancing algorithms{, cf.}~\cite{Biegler.2010}, to tackle time issues, quantum computing potentially presents a paradigm shift, offering hope for efficient real-time solutions even for complex, high-dimensional problems~\cite{symons2023practitioners}.

Quantum computing, and specifically quantum annealing as implemented by D-Wave, has emerged as a promising alternative to classical solvers~\cite{Yarkoni_2022}. 
D-Wave's 5,000-qubit quantum annealer has recently demonstrated a significant quantum advantage in large-scale 3D spin glass optimization problems{~\cite{king2023quantum}. However, it is important to note that this achievement mainly pertains to an academic context, and extending these results to practical, industry-relevant optimization problems is a crucial next step yet to be realized in the realm of quantum computing.}
Quantum annealers, which specialize in optimization and implement physical Ising models, {have the capability to solve NP-Complete Ising instances, and through efficient reformulations, can address a subset of NP-hard optimization problems. 
	As Lucas~\cite{lukas2014quantum} highlighted, the practical applicability to various NP-hard challenges requires mappings that may not be straightforward in all cases}.
Despite their higher qubit count compared to gate-based quantum computers developed by entities like Google and IBM{, it is important to recognize that the operational context and functionality of these qubits are distinct in annealers and gate-based systems. Quantum annealers focus on optimization tasks within the Ising model, whereas gate-based quantum computers offer a more diverse set of quantum logic operations. As such, qubit count is not the sole factor in assessing the computational power or versatility of these different quantum computing paradigms}~\cite{arute2019quantum, kim2023evidence}.

The continuous optimization domain on quantum hardware is still relatively underexplored. 
Binary encoding schemes, commonly utilized in this field, are known to suffer from precision issues on analog Quantum Processing Units (QPUs)~\cite{Yarkoni_2022}. 
{Emerging techniques in quantum computation, such as Quantum BLAS, aim to accelerate linear algebra operations, whereas approaches like the Quantum Approximate Optimization Algorithm (QAOA) with continuous variables are being developed to tackle optimization problems~\mbox{\cite{herman2023quantum, verdon2019quantum, 10.1145/3583133.3596358}}}.
Understanding the challenge landscape, this manuscript explores the application of D-Wave's Advantage Quantum Computer~\cite{king2023quantum} to address dynamic optimization problems. We specifically focus on {a dynamic optimization problem with a} continuously stirred tank reactor (CSTR) {embedded}, presenting a case study as our central investigation.

Dynamic optimization on quantum hardware is almost unexplored.
Deng et al.~\cite{Deng2023} have approached the solution of a dynamic optimization problem for a Heating, Ventilation, and Air Conditioning (HVAC) system process using quantum annealing.
They embedded the analytic solution of the differential equations constituting the process model.
{Fern\'{a}ndez-Villaverde and Hull \cite{fernándezvillaverde2023dynamic} concentrated on a dynamic optimization problem with a real business cycle model, governed by a system of differential-algebraic equations (DAE), embedded. 
	To solve this optimization problem on a quantum computer, they substituted the model DAE with its parametric solution. 
	In effect, their methodology transformed the inherently dynamic optimization problem into a static optimization problem.}

{The existing methodologies to solve dynamic optimization problems on quantum computers~\cite{Deng2023,fernándezvillaverde2023dynamic} are effective in the scarce instances where an analytic or parametric solution of the process model, typically a DAE, is accessible in a form that can be solved on quantum computers.}
{However, they falter when such solutions remain elusive or non-existent.} 
{Generally, deriving analytic or parametric solutions for process models is a formidable challenge, especially when nonlinear equations are involved, e.g., distillation column models}~\cite{Biegler.2010}.
{Furthermore, analytic or parametric solutions often comprise nonlinear terms even in cases of linear differential equations that hinder a treatment on quantum computers, hence preventing a direct embedding those analytic or parametric solutions. This motivates embedding DAEs rather than their solution in the optimization problem to be solved on quantum computers.}

In contrast to the {prior studies~\mbox{\cite{Deng2023,fernándezvillaverde2023dynamic}}}, we embed {a system of differential equations} in the optimization problem and do, in turn, not require the analytic {or parametric} solution of the process model. 
This effort positions us at the frontier of quantum computing applications in the process industry, a sector {that is actively exploring the multifaceted advancements and recognizing the potential challenges and opportunities of quantum computational methods, as outlined by Bernal et al.~\mbox{\cite{bernal2022perspectives}}}.
We navigate through challenges including minor embedding and mapping binary quadratic problems to QPU topology. 
Our analysis focuses on the quantitative performance, specifically the simulation error, of the various methods in solving the optimization problem, as well as on a more qualitative comparison of computational efficiency. {To this end, we employ the gap-at-time metric, which compares the quality of the results obtained after a specific duration of simulation time~\cite{lubinski2024optimization}.}

\subsection*{Outline}
The paper first introduces the class of optimization problems under consideration, followed by a presentation of the CSTR model. We then delve into the classical results, outlining the transformation process into a Quadratic Unconstrained Binary Optimization (QUBO) problem, and sharing experiences encountered with a brute-force solver. Subsequently, we delve into discussions on simulated annealing, quantum annealing, and hybrid solver approaches, accompanied by an evaluation of errors and performance metrics. The paper concludes with a comprehensive assessment and potential directions for future applications of quantum computing in tackling more complex dynamic optimization problems.

\section{Dynamic Optimization Problem Formulation}
\label{sec:optProblem}

{We consider an optimization problem with DAE embedded (OP-DAE).
	The discrete-time OP-DAE} of interest in this study is formally defined as:
\begin{subequations}
	\begin{align}
		\min_{\substack{\boldsymbol{x}_1,...,\boldsymbol{x}_{N} \in \mathcal{X} \\
				\boldsymbol{y}_0,...,\boldsymbol{y}_{N} \in \mathcal{Y}  \\ \boldsymbol{u}_0,...,\boldsymbol{u}_{N} \in \mathcal{U}  } } &  E(\boldsymbol{x},\boldsymbol{y},\boldsymbol{u}), \\
		s.t. \qquad  & \boldsymbol{x}_{t+1} = \boldsymbol{f} (\boldsymbol{x}_t,\boldsymbol{y}_t,\boldsymbol{u}_t), \forall t \in \{ 0,...,N-1 \}, \label{eq:optimizationProblem:dynamics} \\
		& \boldsymbol{0} = \boldsymbol{g} (\boldsymbol{x}_t,\boldsymbol{y}_t,\boldsymbol{u}_t), \forall t \in \{ 0,...,N \}, \label{eq:optimizationProblem:algebraics} \\ & \boldsymbol{x}_0 = \boldsymbol{X_0} \label{eq:optimizationProblem:initials} ,\\
		& \boldsymbol{0} = \boldsymbol{h}^\mathrm{e} (\boldsymbol{x}_t,\boldsymbol{y}_t,\boldsymbol{u}_t) , t\in \{ 0,...,N\}\label{eq:optimizationProblem:equalities} , \\
		& \boldsymbol{0} \geq \boldsymbol{h}^\mathrm{i} (\boldsymbol{x}_t,\boldsymbol{y}_t,\boldsymbol{u}_t) , t\in \{ 0,...,N\}
		\label{eq:optimizationProblem:inequalities},
	\end{align}
	\label{eq:optimizationProblem}
\end{subequations} 
with $\mathcal{X} \subseteq \mathbb{R}^{n_\mathrm{x}}$, $\mathcal{Y} \subseteq \mathbb{R}^{n_\mathrm{y}}$, and $\mathcal{U} \subseteq \mathbb{R}^{n_\mathrm{u}}$, differential states $\boldsymbol{x}$, algebraic states $\boldsymbol{y}$, manipulated variables $\boldsymbol{u}$ number of differential states $n_\mathrm{x}$, number of algebraic states $n_\mathrm{y}$, and number of manipulated variables $n_\mathrm{u}$. 
The functions $E: \mathcal{X} \times \mathcal{Y} \times \mathcal{U} \rightarrow \mathbb{R}$, $\boldsymbol{f}: \mathcal{X} \times \mathcal{Y} \times \mathcal{U} \rightarrow \mathbb{R}^{n_\mathrm{x}}, \boldsymbol{g}: \mathcal{X} \times \mathcal{Y} \times \mathcal{U} \rightarrow \mathbb{R}^{n_\mathrm{y}}, \boldsymbol{h}^e: \mathcal{X} \times \mathcal{Y} \times \mathcal{U} \rightarrow \mathbb{R}^{n_\mathrm{e}}$, and $\boldsymbol{h}^i: \mathcal{X} \times \mathcal{Y} \times \mathcal{U} \rightarrow \mathbb{R}^{n_\mathrm{i}}$ define the objective function, system dynamics, algebraic constraints, and equality and inequality constraints, respectively, with the number of equality constraints $n_\mathrm{e}$ and the number of inequality constraints $n_\mathrm{i}$.
{The inequality constraints comprise $2 \times N^\mathrm{u}$ lower and upper bound constraints for the manipulated variables $u$.}
The initial states are given by $\boldsymbol{X_0} \in \mathcal{X}$.

The dynamic behavior of the system is governed by a DAE with differential index 1, as expressed in Equations~\eqref{eq:optimizationProblem:dynamics} and~\eqref{eq:optimizationProblem:algebraics}. 
When initially presented in continuous-time form, converting the DAEs to a discretized framework is necessary, which can be achieved using strategies outlined in~\cite{Biegler.2010}, such as collocation or Euler methods.
{Notice that we discriminate between the model equations (Equations~\eqref{eq:optimizationProblem:dynamics} and~\eqref{eq:optimizationProblem:algebraics}) and the equality constraints (Equation \eqref{eq:optimizationProblem:equalities}). 
	While all of these equations are formally equality constraints, the model equations have to be satisfied at any time point on the simulation horizon, while the equality constraints (Equation \eqref{eq:optimizationProblem:equalities}) can be evaluated at distinct points in time, e.g., as in the case for endpoint constraints, evaluated at the final time point $t=N$. 
	Furthermore, the equations are handled differently in different solution approaches. 
	For instance, in single shooting, the model equations are (usually) solved using an integrator while the equality constraints are solved by an optimization solver \cite{Biegler.2010}.}

To align the DAE system with the computational capabilities of quantum architectures, we make two key assumptions. First, we assume the system can be converted into an ordinary differential equation (ODE) system, specifically in the form $\boldsymbol{x}_{t+1} = \boldsymbol{f}^\mathrm{ODE}(\boldsymbol{x}_{t}, \boldsymbol{u}_{t})$. 
{This transformation can be achieved by substituting the algebraic variables $\boldsymbol{y}$ in the differential equations by the explicit form of the algebraic equations $\boldsymbol{y}_t = \boldsymbol{g}'(\boldsymbol{x}_{t},\boldsymbol{u}_{t})$ leading to $\boldsymbol{f}^\mathrm{ODE}(\boldsymbol{x}_{t}, \boldsymbol{u}_{t}) = \boldsymbol{f}(\boldsymbol{x}_{t}, \boldsymbol{g}'(\boldsymbol{x}_{t},\boldsymbol{u}_{t}), \boldsymbol{u}_{t})$.}
Second, we require that the function $\boldsymbol{f}^\mathrm{ODE}$ is a polynomial to facilitate its later transformation into a quantum hardware-compatible format. {The direct implementation of non-polynomial functions in binary representations poses considerable challenges for quantum hardware, if not mitigated by approximation strategies~\cite{haner2018optimizing}. In contrast, handling partial differential equations does not present fundamental issues, although the qubit requirements escalate notably as the number of mixed derivatives in the system increases.}

{We} first recast the constrained optimization problem, with ODEs embedded, into an unconstrained one. This is achieved through the application of direct substitution, leading to an unconstrained problem focusing solely on the manipulated variables $u_t$. Subsequent binarization prepares the problem for quantum computational methods.
{To summarize, our approach}
\begin{enumerate}[i]
	\item {transforms an OP-DAE into an optimization problem with ODE embedded (OP-ODE),}
	\item {converts this into an optimization problem with variable bounds only (OP-B) using direct substitution for the variables and penalization for the constraints,}
	\item {binarizes the optimization variables, culminating in a QUBO problem suitable for quantum solutions.}
\end{enumerate}
{The transformation step (i) can easily be done if the algebraic equations of the DAE can be written in an explicit form with respect to the algebraic variables.
	However, this step is a limitation in cases where an explicit form of the algebraic equations is not available.
	Furthermore, step (i) enables faster solution times as the algebraic equations are solved already offline.
	Though the direct substitution step (ii) is more generally applicable than step (i), step (ii) requires that the discrete-time form of the ODE can be reduced by direct substitution leading to an objective function as a function of the initial states and the manipulated variables.
	This is the case, e.g., in the case of an Euler forward or backward discretization.
	Otherwise, the discrete-time ODE has to be added in the form of penalty terms to the objective function for those variables that cannot entirely be substituted.
	Note that step (ii) removes constraints from the optimization problem at the cost of a more complex objective function.
	The binarization and scaling step (iii) requires suitable bounds of the manipulated variables and knowledge about a sufficient discretization of the manipulated variables, as binarization poses a fixed discretization which is sufficient for discrete variables and more challenging for continuous variables.}
The CSTR model, serving as our case study, encapsulates these steps in practice and will be explored in the upcoming section. 
{Table~\ref{tab:problem_overview} summarizes information about the different optimization problems, including number of variables, constraints, and solution times.
}

\begin{table}[t]
	\tymin=80pt
	\footnotesize
	\centering
	\caption{Optimization problem overview. $N^\mathrm{bin}_{u}$ is the number of binary variables used to binarize the manipulated variable $u$. Solution times reported for Ipopt are average values over 10 solution runs. In the chemical reactor example considered in this work, the model is an ODE already. Hence DOP-DAE and DOP-ODE coincide. S. An.: Simulated Annealing. Q. An.: Quantum Annealing.}
	\begin{tabulary}{\textwidth}{L|LLLL} 
		& \textbf{DOP-DAE}  & \textbf{DOP-ODE}  & \textbf{OP-B} & \textbf{QUBO}  \\ \hline
		\textbf{Model \newline embedded} & Original DAE & ODE & no model & no model \\
		& & &  \\
		\textbf{Number of \newline variables} & $N^\mathrm{x}$ differential variables, $N^\mathrm{y}$ algebraic variables, $N^\mathrm{u}$ manipulated variables   & $N^\mathrm{x}$ differential variables, $N^\mathrm{u}$ manipulated variables  & $N^\mathrm{u}$ manipulated variables & $\sum_{u} N^\mathrm{bin}_{u}$ binary variables \\
		\textbf{Number of \newline constraints} & $N^\mathrm{x}$  differential equations, $N^\mathrm{y}$ algebraic equations, $N^\mathrm{e}$ equality constraints, $N^\mathrm{i}$ inequality constraints & $N^\mathrm{x}$  differential equations, $N^\mathrm{e}$ equality constraints, $N^\mathrm{i}$ inequality constraints &  $N^\mathrm{u}$ variable bounds  &  unconstrained problem \\
		\textbf{Optimization \newline solver} & IpOpt & IpOpt & IpOpt & 
		Gurobi, SCIP, MindtPy, S.~An., Q.~An.
		\\
		\textbf{Solution times \newline (wall-clock time)} & 75 ms & 75 ms & 71 ms & see discussion in Sections \ref{sec:simAnnealQubo} to \ref{sec:hybridSolvers} \\ \hline
	\end{tabulary}
	\label{tab:problem_overview}
\end{table}

\section{Optimization with Chemical Reactor Model Embedded}
\label{sec:cstr}

\begin{figure}[t]
	\centering
	\includegraphics[width=0.5\textwidth]{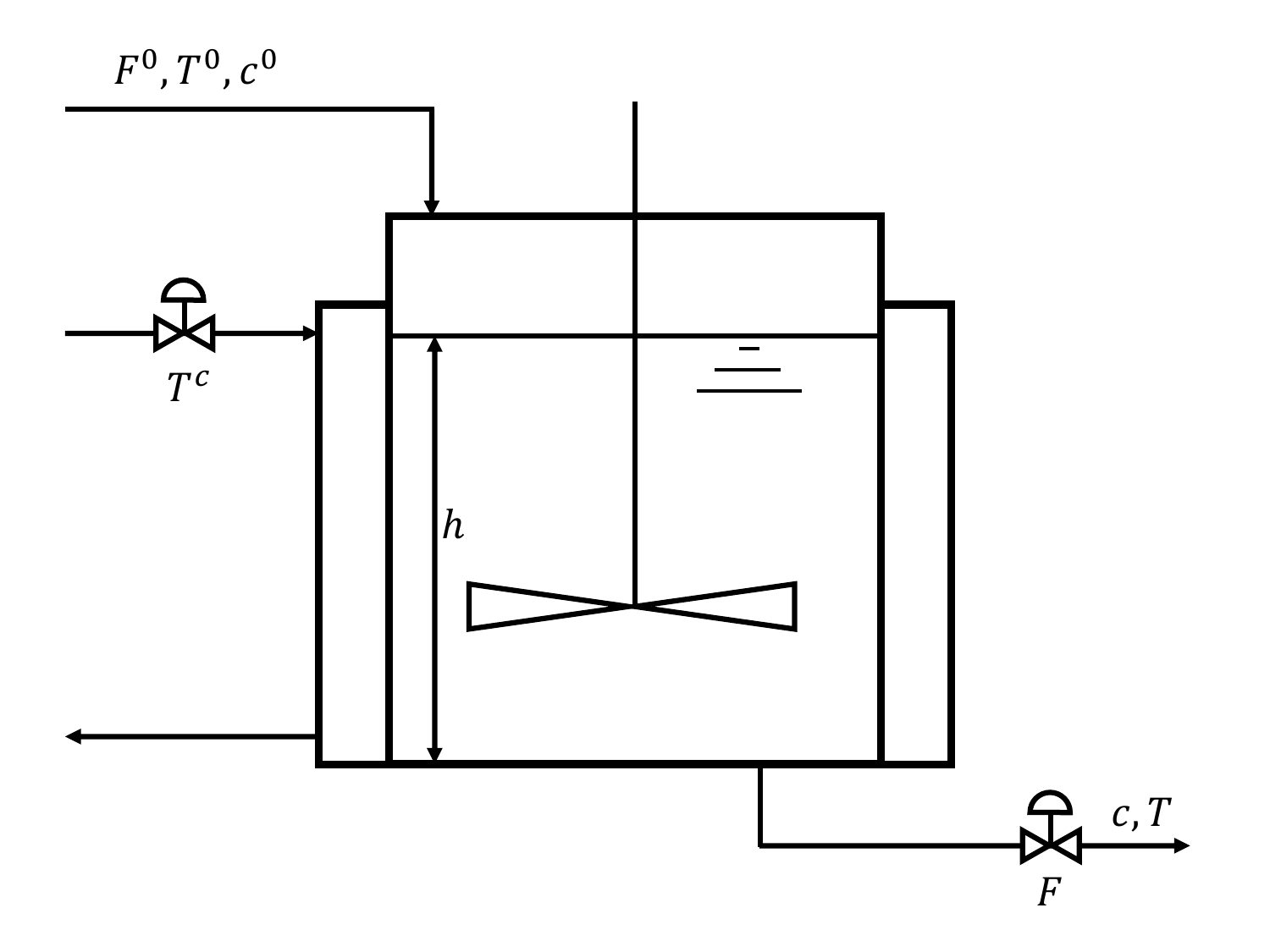}
	\caption{Flowsheet of the Continuous Stirred Tank Reactor with inflow rate $F^0$, temperature $T^0$ and concentration $c^0$. The coolant temperature is denoted by $T^\mathrm{c}$, the height of the reactor by $h$ and the outflow rate by $F$, concentration $c$ and temperature by $T$.}
	\label{fig:reactorFlowsheet}
\end{figure}

The CSTR, a pivotal component in chemical processes, serves as an exemplary model to highlight the challenges and intricacies associated with dynamic optimization in the process industry. 
Its usage as a benchmark for MPC underscores its relevance, making it an ideal candidate for in-depth study and optimization analysis.
We consider the CSTR model as presented in Ref.~\cite{rawlings2017model}.
Figure~\ref{fig:reactorFlowsheet} shows a flowsheet of the CSTR that we use as a case study.
We base our investigation on this specific CSTR model characterized by an Arrhenius term, crucial for capturing the temperature-dependent reaction rates. 
However, considering the nascent stage of quantum computing capabilities, especially in handling and evaluating non-polynomial and exponential terms, we approximate the exponential Arrhenius term as a constant ($T^{\mathrm{fix}} = 340 $K). This simplification is a pragmatic choice, aimed at making the problem more tractable for quantum hardware without expanding the term as a polymial.

The model in its comprehensive form is derived from Reference~\cite{rawlings2017model} and is presented in detail in Appendix~\ref{appendix:chemicalReactorModel}. 
This representation takes the form of a system of a nonlinear continuous-time {ODE}.
{Thus, the DAE-to-ODE transformation of our solution approach can be skipped for the present case study.}
To facilitate the optimization process, we discretize the model using an explicit Euler method, dividing the 4-minute time horizon we use for the dynamic optimization into $N = 20$ equal-spaced time intervals ($\Delta t = 0.2$ min), thereby following a full discretization approach for the solution of the resulting dynamic optimization problem~\cite{Biegler.2010}. 
The following model equations result after the assumptions described in Appendix~\ref{appendix:chemicalReactorModel}:

\begin{subequations}
	\begin{align}
		c_{i+1} & =  c_{i} + \Delta t \cdot  \frac{F^0(c^0 - c_i)}{\pi r^2 h}  \nonumber\\ 
		&\quad - \Delta t \cdot k^0 \exp\left(-\frac{E}{R T^{\mathrm{fix}}}\right) c_{i} \label{eq:simplifiedpannochia:c} \\
		T_{i+1} & = T_i +  \Delta t \cdot \frac{F^0(T^0 - T_i)}{\pi r^2 } \nonumber \\ 
		&\quad -  \Delta t \cdot \frac{\Delta H}{\rho C_\mathrm{p}} k^0 \exp\left(-\frac{E}{R T^{\mathrm{fix}}}\right) c_i \nonumber \\ 
		&\quad  + \Delta t \cdot \frac{2 U}{r \rho C_\mathrm{p}} \left(T^c_i - T_i\right) \label{eq:simplifiedpannochia:T} 
	\end{align}
	\label{eq:simplifiedpannochia}
\end{subequations}

The primary goal is to manipulate the coolant temperature $T^c$ for the reactor temperature $T$ to meet the target value $T^\mathrm{fix}$, directly influencing the reactor's performance, such that we formulate the dynamic optimization problem at hand with the objective function given by
\begin{equation}
	E = \sum_{i=0}^{N}(T_i-T^{\mathrm{fix}})^2 \,,
	\label{eq:energyfunction}
\end{equation}
where the summation extends over all discrete time steps.
The resulting dynamic optimization problem is presented in Appendix~\ref{appendix:DynamicOptimizationProblem}.

\begin{figure}[t]
	\centering
	\includegraphics[width=\textwidth]{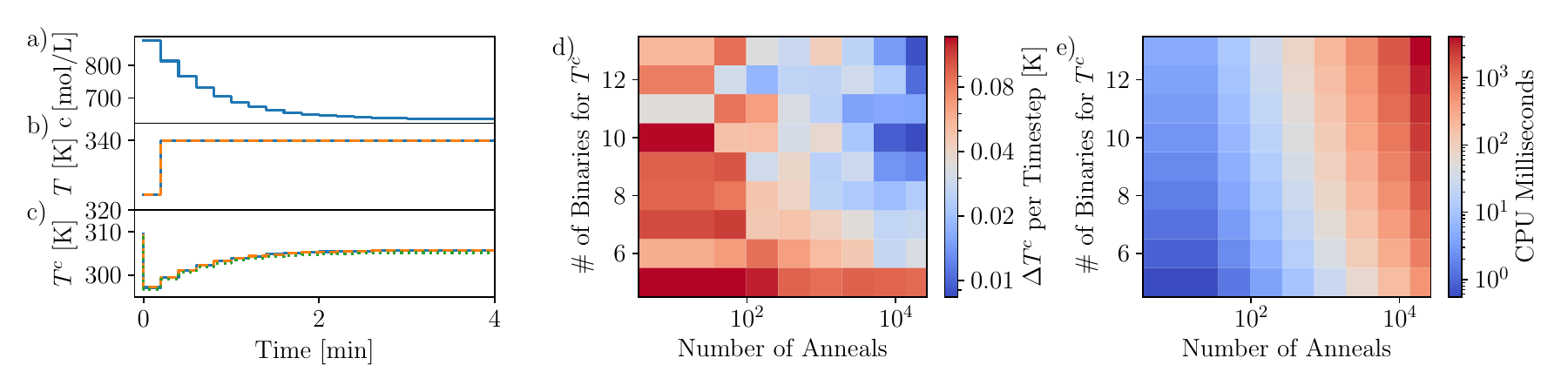}
	\caption{Comparison of concentration (a), reactant temperature (b), and cooling temperature (c) across three distinct methodologies for the reactor case study: the classical Ipopt solution of the DOP-ODE (blue, solid line), the best simulated annealing result for the QUBO (orange, dashed line), and {the Gurobi solution of the QUBO} (green, dotted line). Subplot (d) showcases the optimization error in $T^c$ per timestep, visualized as a colormap, under varying conditions of simulated annealing events and binarization bits for $T^c$. Subplot (e) depicts the computational time necessary to achieve these results, providing a comprehensive performance comparison.}
	\label{fig:simulated_annealing}
\end{figure}

To establish a baseline and for comparative purposes, we solve the optimization problem with DAE embedded, i.e., the original optimization problem of interest, using a classical full-discretization approach~\cite{Biegler.2010}.
Applying an explicit Euler method yields exactly the optimization presented in Appendix~\ref{appendix:DynamicOptimizationProblem}.
To solve this optimization problem, we leverage the Pyomo Dae framework~\cite{Nicholson2017} \footnote{Pyomo Dae, \url{www.pyomo.org}, accessed: 11-08-2023.} with the widely-used solver nonlinear programming solver Ipopt version 3.11.1~\cite{Waechter2005} and solve the resulting optimization problem on classical hardware \footnote{Windows 10 desktop computer with an Intel Core i5-7200U CPU running at 2.50 GHz and equipped with 8 GB RAM.}, with which we reach a converged result (tolerances $10^{-1}$ to $10^{-10}$), i.e., a local optimum within approx. {70 - 80} \,ms {(cf. Table~\ref{tab:problem_overview})}.
The benchmark solution is depicted in Figure~\ref{fig:simulated_annealing}~(a).

Solving the dynamic optimization problem on quantum hardware requires a transformation, since usual quantum hardware cannot treat constrained optimization problems and requires, in turn, a transformation to an unconstrained problem.
The strategy for this transformation entails iteratively substituting all occurrences of $T$ and $c$ in Equation~\eqref{eq:energyfunction}, ultimately transforming the problem into an optimization task exclusively for the coolant temperature $T^c$. 
Although handling the differential equations as actual constraints is a viable alternative, it introduces the complexity of treating $T$ and $c$ as free parameters, complicating the binarization process necessary for quantum hardware, as shown later.
This simplification, while beneficial for our current analysis, highlights a broader challenge in quantum optimization, especially for more complex problems, and sets the stage for further exploration in subsequent sections of this paper.
This approach streamlines the problem for quantum hardware applications in subsequent steps.

\section{Simulated Annealing for QUBO Formulation}
\label{sec:simAnnealQubo}

In this section, we apply simulated annealing as a benchmark to our quantum computing approach for solving the unconstrained optimization problem resulting from the transformations presented in Section~\ref{sec:cstr}.
Using simulated annealing is well-suited as a benchmark since it is a strategy similar to the solution approach on quantum hardware, but running on classical hardware instead of a quantum computer.

Simulated annealing is a probabilistic optimization algorithm originally inspired by the annealing process in metallurgy~\cite{doi:10.1126/science.220.4598.671}. 
It has gained popularity for its simplicity and capability to find near-optimal solutions in large solution spaces with little prior knowledge about the problem. 
However, it is more suited to unconstrained optimization problems and does not guarantee convergence to a global or local optimum.
Simulated annealing works by iteratively improving the solution, allowing for occasional steps that worsen the objective function to escape local minima. 
The probability of accepting such steps decreases over time, mimicking the annealing process as, e.g., implemented in D-Wave's quantum processor.

To prepare for comparison with quantum computing results, we need to transform the unconstrained optimization problem (variables substituted as described in Section \ref{sec:cstr}) into a QUBO, as quantum annealers are designed to solve problems expressed in this format. 
This involves representing our optimization variables, which are continuous, with binary variables.

\begin{figure*}[t]
	\centering
	\includegraphics[width=0.67\textwidth]{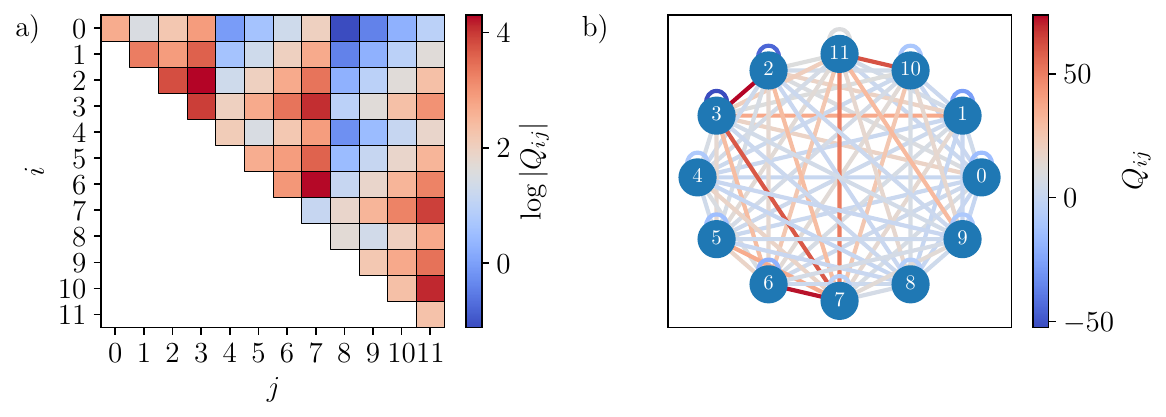}
	\caption{Two representations of the optimization problem for a scenario with three timesteps and a binarization resolution of four bits for the cooling temperature $T^c$: (a) showcases the $Q$ matrix, encapsulating the linear and quadratic coefficients of the problem, while (b) illustrates the corresponding graph network, highlighting the interconnections between binary variables.}
	\label{fig:qubo_matrix}
\end{figure*}

We transform the continuous variables of our problem into binary variables through a discretization process. 
Given a continuous variable \( T_j \) at time-step $j$, we map it to \( n \) binary variables \( b_0, b_1, \ldots, b_{n-1} \), where \( n \) represents the number of bits used for the discretization. This mapping can be expressed as follows:
\begin{equation}
	T_j = T_{\text{lower}} + \left(T_{\text{higher}} - T_{\text{lower}}\right) \cdot \frac{\sum_{i=0}^{n-1} b_i \cdot 2^i}{2^n - 1} \, ,
	\label{eq:binarization}
\end{equation}
where \( T_{\text{lower}} \) and \( T_{\text{higher}} \) correspond to the bounds of the optimization problem, ensuring that all solutions lie within this interval. 
The resolution of this discretization is given by $(T_{\text{higher}} - T_{\text{lower}})/(2^n - 1)$.

Substituting Equation~\eqref{eq:binarization} into the unconstrained optimization problem, we obtain the QUBO formulation of our optimization problem. 
This can be expressed as:
\begin{equation*}
	E = \sum_i Q_{ii} b_i + \sum_{i< j}Q_{ij} b_i b_j,
\end{equation*}
where $Q_{ii}$ represents the linear coefficients, and $Q_{ij}$ (for $i < j$) are the quadratic coefficients. 
The QUBO problem can thus be formulated as:
\begin{equation}
	\min_{\boldsymbol{b}\in \{0,1\}^n} \boldsymbol{b}^T \boldsymbol{Q} \boldsymbol{b} \,.
	\label{eq:qubo_opt}
\end{equation}
Again, achieving this form is possible even when starting with systems of ODEs with higher-order polynomial terms~\cite{dattani2019quadratization}.

To have a further suitable benchmark, we solve the resulting QUBO~\eqref{eq:qubo_opt} using the optimization solver{s} Gurobi version 9.1.0~\footnote{Gurobi Optimization, LLC., \url{www.gurobi.com}, accessed: 11-08-2023.}, {SCIP Optimization Suite 8.0}~\footnote{{Zuse Institute Berlin (ZIB), SCIP Optimization Suite, \url{www.scipopt.org}, accessed: 12-02-2024.}}~{\mbox{\cite{10.1145/3585516}}}, {and MindtPy}~\footnote{{MindPy, \url{https://pyomo.readthedocs.io/en/stable/contributed_packages/mindtpy.html}, accessed: 12-02-2024.}}~{\mbox{\cite{BERNAL2018895}} (with the outer-approximation algorithm), all} interfaced via Pyomo ~\cite{bynum2020pyomo, hart2011pyomo} \footnote{Pyomo, \url{www.pyomo.org}, accessed: 11-08-2023.}.

\begin{figure*}[t]
	\centering
	\includegraphics[width=\textwidth]{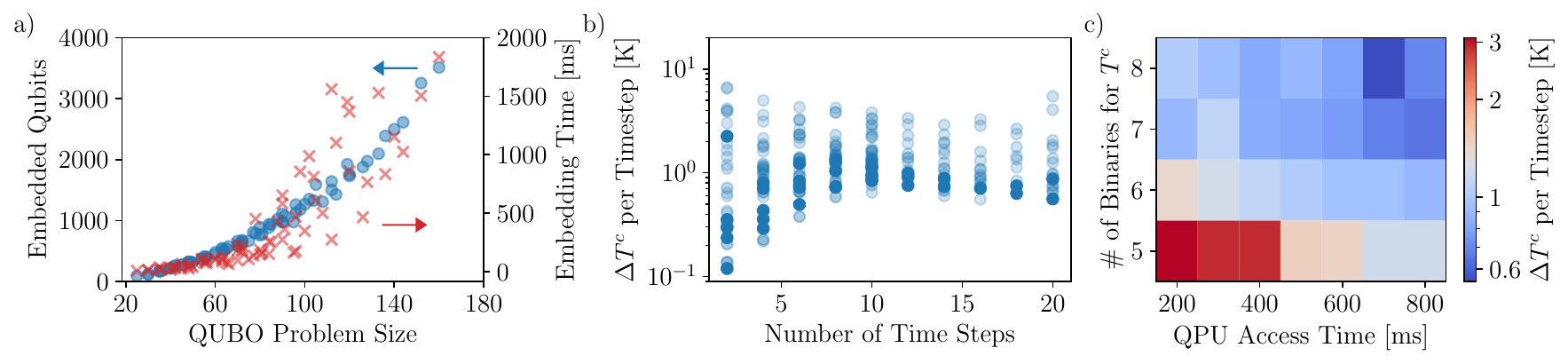}
	\caption{Results for the CSTR model using D-Wave’s Advantage system. (a) Illustrates the number of embedded qubits {(blue dots)} as a function of the original QUBO problem size {(number of binary variables)}, which showcases the increased complexity and qubit requirements for larger problems. {Additionally, it shows the embedding time on classical hardware (red crosses).} (b) Depicts the error in $T^c$ per timestep across various problem sizes {(number of binaries per $T^c$, number of anneals, and anneal time, as well as different chain strengths between 50 and 5000)}, where the color intensity corresponds to the QPU access time; longer compute-times (darker) generally correlate with lower errors. (c) Provides a colormap representation of the error in $T^c$ per timestep for a set of fixed simulation parameters {(chain strength fixed at 50, annealing time fixed at 15 $\mu$s) and varying the number of anneals to achieve different total QPU access times, which are grouped into the illustrated timeslots,} demonstrating a trend towards improved results for longer computations and better representations (number of binaries for $T^c$).}
	\label{fig:real_qc_results}
\end{figure*}

The matrix $\boldsymbol{Q}$ in the QUBO~\eqref{eq:qubo_opt} , exemplified for a smaller instance of the problem, is visualized in Figure~\ref{fig:qubo_matrix}.

In this format, the optimization problem is now suitable for solution using quantum annealers. 
It is worth mentioning that if we encountered terms of higher order beyond quadratic, it is often possible to convert them into quadratic terms at the cost of introducing auxiliary variables~\cite{dattani2019quadratization}. 

The results of applying simulated annealing, {as implemented in D-Wave's Ocean software}~\footnote{D-Wave Samplers v. 1.2.0, \url{https://docs.ocean.dwavesys.com/en/latest/docs_samplers/release_notes.html}, accessed: 02-02-2024}, to our QUBO problem~\eqref{eq:qubo_opt} are presented in Figure~\ref{fig:simulated_annealing}. Figure~\ref{fig:simulated_annealing}~(c) shows the optimal cooling temperature $T^c$ obtained through different methodologies: the classical solution of the original problem with DAE embedded using Pyomo Dae with IPOPT, the optimal simulated annealing result for the QUBO~\eqref{eq:qubo_opt}, and the Gurobi {solution} for the QUBO~\eqref{eq:qubo_opt}.  
Qualitatively, all three results overlap for the entire simulation period. 
A deeper analysis for the full problem with 200 binary variables {(using ten binary variables to represent the cooling temperature at each of the 20 timesteps)}, shown in Figure~\ref{fig:simulated_annealing}(d), reveals that the optimization error in $T^c$ per timestep is within the permyriad range for the best results, increasing for lower numbers of bits used for binarization and lower numbers of annealing steps, demonstrating that the simulated annealing approach for the QUBO formulation can yield results comparable to traditional optimization techniques, given sufficient computational resources. We maintain this metric for error assessment throughout the manuscript, as it offers a tangible gauge of result quality. However, as depicted in panel (e), the computational time required increases drastically with both the number of bits for binarization and the number of annealing steps, highlighting the trade-off between solution quality and computational efficiency. 
Compared to the classical solution of the original dynamic optimization problem {(OP-DAE/ODE)} using Ipopt{, we achieve a deviation of approx. 30 mK (8 mK) per timestep below 100 ms (1 s) compute time}.
Note that {the solvers} {SCIP and MindtPy} did not provide a converged and useful solution to the QUBO problem within the time limit of {6} h, indicating that the computational challenge originates nature of the QUBO problem statement. {On the other hand, Gurobi solved the QUBO in about 610 ms wall-clock time (time averaged over 10 solution runs), setting the benchmark for our exploration on real quantum hardware.}

\section{Quantum Annealing on D-Wave's Advantage Quantum Computing System}
To solve a QUBO directly on a quantum annealer, a mapping called minor embedding is required to fit the problem onto the quantum processing unit’s (QPU) topology~\cite{zbinden2020embedding}. 
In D-Wave's Advantage Pegasus topology~\cite{dattani2019pegasus}, most qubits are connected to fifteen other qubits, allowing for complex problem mappings. 
Some algorithms can embed a problem of \(N\) variables in at most \(N^2\) qubits, although the efficiency of the embedding can significantly impact performance. 
In this study, the embedding was found using a heuristic approach for minor embedding as presented by Cai et al.~\cite{cai2014practical}.

Important parameters for our study include, but are not limited to, the annealing time, the number of annealing runs, and the chain strength. 
The chain strength, in particular, determines the effectiveness of the embedding, ensuring that when it is necessary to represent a logical qubit by multiple physical qubits (linked by a chain), e.g. when the problem connectivity exceeds the physical qubit connections, they assume the same value at the end of the annealing cycle~\footnote{D-Wave Systems Inc., ``Programming the D-Wave QPU: Setting the Chain Strength'', \url{www.dwavesys.com/media/vsufwv1d/14-1041a-a_setting_the_chain_strength.pdf}, accessed: 09-30-2023. \label{fn:dwave_chain}}.  
If the chain strength is not adequately chosen, it can lead to broken chains and suboptimal solutions, necessitating postprocessing or re-runs with adjusted parameters. 

In our quantum annealing experiments {performed on D-Wave's Advantage 6.4 system}, we analyzed scenarios with {up to} eight binary variables per discrete cooling temperature $T^c$. The outcomes, as depicted in Figure~\ref{fig:real_qc_results}, provide a comparative analysis across different parameters and settings.
A significant challenge highlighted in these results is the rapid escalation in the number of required qubits for larger problem sizes due to the embedding process. As shown in Figure~\ref{fig:real_qc_results}(a){the largest problems we were able to fit onto the QPU consist of eight binary variables per timestep; that is, 160 variables describe the full-time horizon. Beyond this problem size, the embedding algorithm was unable to find an embedding on the solver topology. Additionally, the time needed to find the embedding on conventional hardware quickly exceeds the calculation times accessible to us on quantum hardware (approximately 1000 ms), especially for the larger problems, adding to the total solution time}. The underlying {quadratically bounded} growth in qubit requirements underscores the computational hurdles and the need for ongoing research in efficiently mapping problems onto quantum hardware. {We visualize the embedding for two examples in Appendix~\ref{appendix:pegasusembedding}·}

{Despite the intricate settings and the capabilities of the quantum hardware, our results also indicate that the error in solutions with respect to solutions obtained by solving the DOP-DAE, DOP-ODE, and OP-B using Ipopt and the QUBO using Gurobi is generally larger than those obtained using the simulated annealing results.}
This is the case even under optimal conditions, where errors reach values within the percent range of the coolant temperature. Notably, smaller problem sizes, which use fewer timesteps, typically allow for longer simulations, leading to reduced errors as detailed in Figure~\ref{fig:real_qc_results}(b). 

{To obtain the results shown in Figure~\ref{fig:real_qc_results}(c), we varied the number of annealing events (ranging from 1000 to 4000) and the number of binary variables representing $T^c$, while holding the chain strength and annealing time constant at 50 and 15 $\mu$s, respectively, and maintaining all other settings at their default values. We observed no significant improvement in our results with longer annealing times. QPU access time refers to the total duration allocated on the D-Wave Advantage chip for preparing, executing, and reading out one simulation run}~\footnote{D-Wave Systems Inc., ``Operation and Timing'', \url{https://docs.dwavesys.com/docs/latest/c_qpu_timing.html}, accessed: 02-02-2024. \label{fn:dwave_timing}}. {While we achieved errors as small as 0.55 K per timestep in one of the longest simulations using eight binary variables per timestep, this required the maximum time permitted for public access to the D-Wave chip. Furthermore, the largest problem size of 160 qubits necessitated nearly 2 seconds for embedding alone, rendering a speed comparison between quantum and classical solvers irrelevant when considering both embedding time and QPU access time.}

\section{Hybrid Solvers}
\label{sec:hybridSolvers}
\begin{figure*}[t]
	\centering
	\includegraphics[width=\textwidth]{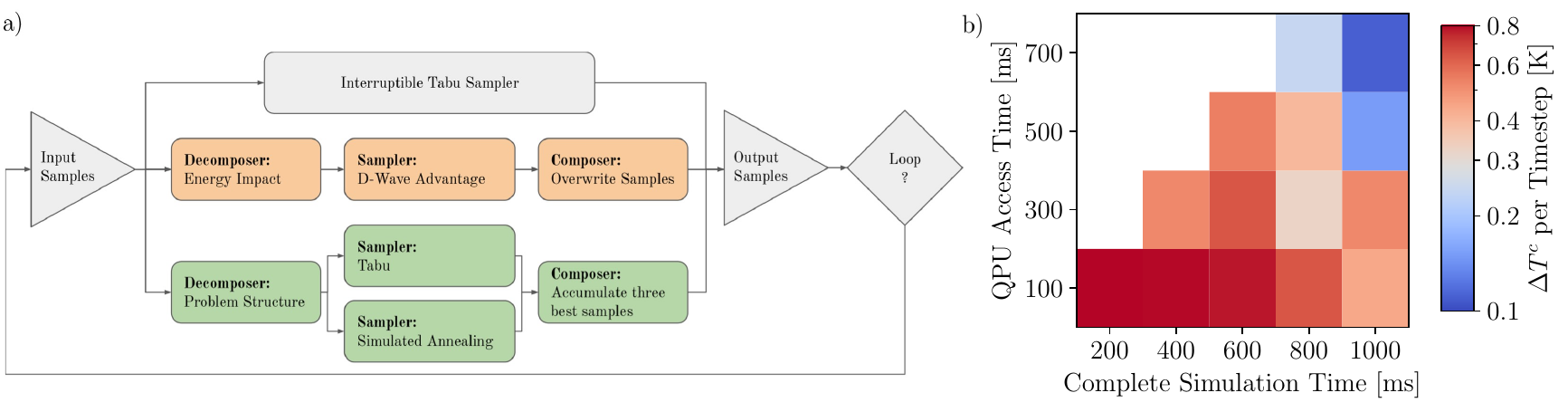}
	\caption{(a) displays a flowchart following Footnote~\ref{fn:dwave_hybrid}, illustrating the Kerberos algorithm's iterative routine. An initial input undergoes simultaneous computation through various solvers: a classical interruptible tabu search addressing the full problem, running in tandem with branches that subdivide the problem employing solvers on quantum processors and additional setups of parallel simulated annealing and tabu search. Decomposers fragment the problem, samplers tackle each piece, and composers merge the partial solutions back into a whole. This concurrent activity yields an optimized sample collection for the next iteration. (b) A colormap demonstrating the error in coolant temperature $T^c$ across each time step in relation to the count of anneals and annealing times per cycle, revealing improved performance with more extended annealing durations. {Blank fields indicate a lack of data due to incompatibility between QPU access time and total runtime.}}
	\label{fig:hybrid_qc_results}
\end{figure*}

{To mitigate the long embedding times associated with large problem sizes, we employed a hybrid approach designed to divide a larger problem into several smaller subproblems.}  Utilizing D-Wave's Kerberos, an asynchronous hybrid decomposition sampler, we address the optimization problem by harnessing quantum acceleration alongside the capability to process problems of any scale~\footnote{D-Wave Systems Inc., ``D-wave hybrid'',  \url{https://docs.ocean.dwavesys.com/projects/hybrid/en/latest}, accessed: 11-08-2023. \label{fn:dwave_hybrid}}. 
This method presents a promising direction to potentially surpass classical solvers in certain instances~\cite{PhysRevA.105.062435,stogiannos2022experimental}, negating the necessity for extensive pre-computation of embeddings. The Kerberos hybrid sampler operates by executing three parallel sampling branches {(two classical sampling routines and one quantum solver)} and its procedure is detailed in Figure~\ref{fig:hybrid_qc_results}(a).

In our experiments for the full problem with 160 {binary variables}, we varied the number of annealing runs, {the number of iterations through the loop in Figure~\ref{fig:hybrid_qc_results}(a), the time limit for both tabu sampler and simulated annealing} and limiting the subproblem size to 40 variables. This calibration aimed to synchronize the operating intervals between the tabu sampler and the quantum processor, ensuring QPU {access times between 100 ms and 700 ms through varying the number of annealing steps and overall simulation times between 200 ms and 1000 ms}. The outcomes, depicted in Figure~\ref{fig:hybrid_qc_results}(b), illustrate the hybrid quantum solver's superiority over standalone quantum annealing concerning accuracy, yielding more precise results within comparable timeframes. This vindicates the hybrid solver's efficacy in dynamic optimization tasks, integrating the advantages of both classical and quantum computing techniques. {A comparison with simulated annealing results of comparable runtime, however, shows lower performance, which attribute to the splitting of the problem into smaller-sized subproblems and the overhead generated through this process.}



\section{Summary and Future Prospects}
Our exploration of dynamic optimization problems utilizing a simplified chemical reactor model has shed light on the formidable challenges involved in transitioning quantum computing from theoretical constructs to practical utility. 
{There are dynamic optimization problems challenging to solve even using state-of-the-art approaches on conventional computers, e.g., \cite{Caspari2019f, Caspari2019b}. 
	Here, we focused on a simple example to show which steps are generally required to solve dynamic optimization problems on quantum hardware and highlight the current opportunities and challenges in this field.}
Although we made substantial simplifications, the scaling properties of the resulting QUBO format on quantum hardware remain unfavorable, underscoring the inherent complexity of dynamic optimization.

Nevertheless, experimental implementations on real quantum hardware have proven that it is feasible to achieve practical outcomes with errors in the percent range within confined time frames. 
Hybrid Quantum Annealing techniques have demonstrated their capacity for managing problems that manage the challenges of embedding, a pivotal factor for applicability in industrial contexts.

The evolution of quantum hardware{, particularly in terms of qubit number and connectivity,} is pivotal for the efficient resolution of larger and more intricate dynamic optimization problems. Critically, optimization algorithms that inherently accommodate constraints without resorting to penalty terms in the objective function or explicit variable substitutions optimized for QUBO representations are highly desirable. 
Moreover, transforming dynamic optimization problems, akin to the model we assessed, to an annealing architecture with a limited number of connections results in a proliferation of auxiliary qubits necessary for the mapping, which impairs the computational process.

Optimizing our methodology using refined discretization strategies that capitalize on continuous variable correlations represents an avenue with considerable potential~\cite{doi:10.7566/JPSJ.92.044802}. Further investigation is imperative to decipher the intricacies of embedding complex dynamic optimization problems, such as DAEs with general nonlinear functions, on quantum hardware, and to examine the nature of dynamic optimization problems with binary variables{, e.g. optimal routing for hybrid electric vehicles~\cite{Caspari_routing_2022}}, which introduce a higher level of complexity for classical approaches.

\section*{Acknowledgments}
The authors thank Dr. Manuel Dahmen from Forschungszentrum Jülich, Germany and Prof. Alexander Mitsos, Ph.D. from RWTH Aachen University, Germany for fruitful discussions.

\bibliographystyle{ieeetr}

\begin{thebibliography}{10}
	
	\bibitem{Biegler.2010}
	L.~T. Biegler, {\em Nonlinear Programming: Concepts, Algorithms, and
		Applications to Chemical Processes}.
	\newblock SIAM, 2010.
	
	\bibitem{rawlings2017model}
	J.~B. Rawlings, D.~Q. Mayne, and M.~Diehl, {\em Model predictive control:
		theory, computation, and design}, vol.~2.
	\newblock Nob Hill Publishing Madison, WI, 2017.
	
	\bibitem{Caspari2019b}
	A.~Caspari, C.~Offermanns, P.~Schäfer, A.~Mhamdi, and A.~Mitsos, ``A flexible
	air separation process: 2. optimal operation using economic model predictive
	control,'' {\em {AIChE} Journal}, vol.~65, no.~11, 2019.
	
	\bibitem{Caspari2019g}
	A.~Caspari, C.~Tsay, A.~Mhamdi, M.~Baldea, and A.~Mitsos, ``The integration of
	scheduling and control: Top-down vs. bottom-up,'' {\em Journal of Process
		Control}, vol.~91, pp.~50 -- 62, 2020.
	
	\bibitem{Schaefer2019a}
	P.~Schäfer, A.~Caspari, A.~Mhamdi, and A.~Mitsos, ``Economic nonlinear model
	predictive control using hybrid mechanistic data-driven models for optimal
	operation in real-time electricity markets: In-silico application to air
	separation processes,'' {\em Journal of Process Control}, vol.~84,
	pp.~171--181, 2019.
	
	\bibitem{Schulze2020}
	J.~C. Schulze, A.~Caspari, C.~Offermanns, A.~Mhamdi, and A.~Mitsos, ``Nonlinear
	model predictive control of ultra-high-purity air separation units using
	transient wave propagation model,'' {\em Computers {\&} Chemical
		Engineering}, p.~107163, 2020.
	
	\bibitem{Kadam2004}
	J.~V. Kadam and W.~Marquardt, ``Sensitivity-based solution updates in
	closed-loop dynamic optimization,'' {\em {IFAC} Proceedings Volumes},
	vol.~37, no.~9, pp.~947--952, 2004.
	
	\bibitem{Vaupel2020}
	Y.~Vaupel, N.~C. Hamacher, A.~Caspari, A.~Mhamdi, I.~G. Kevrekidis, and
	A.~Mitsos, ``Accelerating nonlinear model predictive control through machine
	learning,'' {\em Journal of Process Control}, vol.~92, pp.~261 -- 270, 2020.
	
	\bibitem{symons2023practitioners}
	B.~C.~B. Symons, D.~Galvin, E.~Sahin, V.~Alexandrov, and S.~Mensa, ``A
	practitioner's guide to quantum algorithms for optimisation problems,'' {\em
		arXiv preprint arXiv:2305.07323}, 2023.
	
	
	\bibitem{king2023quantum}
	A.~D. King, J.~Raymond, T.~Lanting, R.~Harris, A.~Zucca, F.~Altomare, A.~J.
	Berkley, K.~Boothby, S.~Ejtemaee, C.~Enderud, {\em et~al.}, ``Quantum
	critical dynamics in a 5,000-qubit programmable spin glass,'' {\em Nature},
	vol.~617, pp.~61--66, 2023.
	
	\bibitem{lukas2014quantum}
	A.~Lucas, ``Ising formulations of many np problems,'' {\em Frontiers in
		Physics}, vol.~2, 2014.
	
	\bibitem{arute2019quantum}
	F.~Arute, K.~Arya, R.~Babbush, D.~Bacon, J.~C. Bardin, R.~Barends, R.~Biswas,
	S.~Boixo, F.~G. Brandao, D.~A. Buell, {\em et~al.}, ``Quantum supremacy using
	a programmable superconducting processor,'' {\em Nature}, vol.~574, no.~7779,
	pp.~505--510, 2019.
	
	\bibitem{kim2023evidence}
	Y.~Kim, A.~Eddins, S.~Anand, K.~X. Wei, E.~Van Den~Berg, S.~Rosenblatt,
	H.~Nayfeh, Y.~Wu, M.~Zaletel, K.~Temme, {\em et~al.}, ``Evidence for the
	utility of quantum computing before fault tolerance,'' {\em Nature},
	vol.~618, no.~7965, pp.~500--505, 2023.
	
	\bibitem{Yarkoni_2022}
	S.~Yarkoni, E.~Raponi, T.~Bäck, and S.~Schmitt, ``Quantum annealing for
	industry applications: introduction and review,'' {\em Reports on Progress in
		Physics}, vol.~85, p.~104001, 2022.
	
	\bibitem{herman2023quantum}
	D.~Herman, C.~Googin, X.~Liu, Y.~Sun, A.~Galda, I.~Safro, M.~Pistoia, and
	Y.~Alexeev, ``Quantum computing for finance,'' {\em Nature Reviews Physics},
	pp.~1--16, 2023.
	
	\bibitem{verdon2019quantum}
	G.~Verdon, J.~M. Arrazola, K.~Brádler, and N.~Killoran, ``A quantum
	approximate optimization algorithm for continuous problems,'' {\em arXiv
		preprint arXiv:1902.00409}, 2019.
	
	\bibitem{10.1145/3583133.3596358}
	J.~Stein, F.~Chamanian, M.~Zorn, J.~N\"{u}\ss{}lein, S.~Zielinski,
	M.~K\"{o}lle, and C.~Linnhoff-Popien, ``Evidence that pubo outperforms qubo
	when solving continuous optimization problems with the qaoa,'' in {\em
		Proceedings of the Companion Conference on Genetic and Evolutionary
		Computation}, GECCO '23 Companion, (New York, NY, USA), p.~2254–2262,
	Association for Computing Machinery, 2023.
	
	\bibitem{Deng2023}
	Z.~Deng, X.~Wang, and B.~Dong, ``Quantum computing for future real-time
	building hvac controls,'' {\em Applied Energy}, vol.~334, p.~120621, 2023.
	
	\bibitem{fernándezvillaverde2023dynamic}
	J.~Fernández-Villaverde and I.~Hull, ``Dynamic programming on a quantum
	annealer: Solving the {RBC} model,'' {\em arXiv preprint arXiv:2306.04285},
	2023.
	
	\bibitem{bernal2022perspectives}
	D.~E. Bernal, A.~Ajagekar, S.~M. Harwood, S.~T. Stober, D.~Trenev, and F.~You,
	``Perspectives of quantum computing for chemical engineering,'' {\em AIChE
		Journal}, vol.~68, no.~6, p.~e17651, 2022.
	
	\bibitem{lubinski2024optimization}
	T.~Lubinski, C.~Coffrin, C.~McGeoch, P.~Sathe, J.~Apanavicius, and D.~E.~B.
	Neira, ``Optimization applications as quantum performance benchmarks,'' 2024.
	
	\bibitem{haner2018optimizing}
	T.~H{\"a}ner, M.~Roetteler, and K.~M. Svore, ``Optimizing quantum circuits for
	arithmetic,'' {\em arXiv preprint arXiv:1805.12445}, 2018.
	
	\bibitem{Nicholson2017}
	B.~Nicholson, J.~D. Siirola, J.-P. Watson, V.~M. Zavala, and L.~T. Biegler,
	``pyomo.dae: a modeling and automatic discretization framework for
	optimization with differential and algebraic equations,'' {\em Mathematical
		Programming Computation}, vol.~10, no.~2, pp.~187--223, 2017.
	
	\bibitem{Waechter2005}
	A.~Wächter and L.~T. Biegler, ``On the implementation of an interior-point
	filter line-search algorithm for large-scale nonlinear programming,'' {\em
		Mathematical Programming}, vol.~106, no.~1, pp.~25--57, 2005.
	
	\bibitem{doi:10.1126/science.220.4598.671}
	S.~Kirkpatrick, C.~D. Gelatt, and M.~P. Vecchi, ``Optimization by simulated
	annealing,'' {\em Science}, vol.~220, no.~4598, pp.~671--680, 1983.
	
	\bibitem{dattani2019quadratization}
	N.~Dattani, ``Quadratization in discrete optimization and quantum mechanics,''
	{\em arXiv preprint arXiv:1901.04405}, 2019.
	
	\bibitem{10.1145/3585516}
	K.~Bestuzheva, M.~Besan\c{c}on, W.-K. Chen, A.~Chmiela, T.~Donkiewicz, J.~van
	Doornmalen, L.~Eifler, O.~Gaul, G.~Gamrath, A.~Gleixner, L.~Gottwald,
	C.~Graczyk, K.~Halbig, A.~Hoen, C.~Hojny, R.~van~der Hulst, T.~Koch,
	M.~L\"{u}bbecke, S.~J. Maher, F.~Matter, E.~M\"{u}hmer, B.~M\"{u}ller, M.~E.
	Pfetsch, D.~Rehfeldt, S.~Schlein, F.~Schl\"{o}sser, F.~Serrano, Y.~Shinano,
	B.~Sofranac, M.~Turner, S.~Vigerske, F.~Wegscheider, P.~Wellner, D.~Weninger,
	and J.~Witzig, ``Enabling research through the scip optimization suite 8.0,''
	{\em ACM Trans. Math. Softw.}, vol.~49, 2023.
	
	\bibitem{BERNAL2018895}
	D.~E. Bernal, Q.~Chen, F.~Gong, and I.~E. Grossmann, ``Mixed-integer nonlinear
	decomposition toolbox for pyomo (mindtpy),'' in {\em 13th International
		Symposium on Process Systems Engineering (PSE 2018)} (M.~R. Eden, M.~G.
	Ierapetritou, and G.~P. Towler, eds.), vol.~44 of {\em Computer Aided
		Chemical Engineering}, pp.~895--900, Elsevier, 2018.
	
	\bibitem{bynum2020pyomo}
	M.~L. Bynum, G.~Hackebeil, W.~E. Hart, C.~D. Laird, B.~L. Nicholson, J.~D.
	Siirola, J.-P. Watson, and D.~L. Woodruff, ``Pyomo-optimization modeling in
	python 3rd ed.,'' tech. rep., Sandia National Lab.(SNL-NM), Albuquerque, NM
	(United States), 2020.
	
	\bibitem{hart2011pyomo}
	W.~E. Hart, J.-P. Watson, and D.~L. Woodruff, ``Pyomo: modeling and solving
	mathematical programs in python,'' {\em Mathematical Programming
		Computation}, vol.~3, pp.~219--260, 2011.
	
	\bibitem{zbinden2020embedding}
	S.~Zbinden, A.~B{\"a}rtschi, H.~Djidjev, and S.~Eidenbenz, ``Embedding
	algorithms for quantum annealers with chimera and pegasus connection
	topologies,'' in {\em International Conference on High Performance
		Computing}, pp.~187--206, Springer, 2020.
	
	\bibitem{dattani2019pegasus}
	N.~Dattani, S.~Szalay, and N.~Chancellor, ``Pegasus: The second connectivity
	graph for large-scale quantum annealing hardware,'' {\em arXiv preprint
		arXiv:1901.07636}, 2019.
	
	\bibitem{cai2014practical}
	J.~Cai, W.~G. Macready, and A.~Roy, ``A practical heuristic for finding graph
	minors,'' {\em arXiv preprint arXiv:1406.2741}, 2014.
	
	\bibitem{PhysRevA.105.062435}
	Y.~W. Koh and H.~Nishimori, ``Quantum and classical annealing in a continuous
	space with multiple local minima,'' {\em Phys. Rev. A}, vol.~105, p.~062435, 2022.
	
	\bibitem{stogiannos2022experimental}
	E.~Stogiannos, C.~Papalitsas, and T.~Andronikos, ``Experimental analysis of
	quantum annealers and hybrid solvers using benchmark optimization problems,''
	{\em Mathematics}, vol.~10, no.~8, p.~1294, 2022.
	
	\bibitem{Caspari2019f}
	A.~Caspari, S.~R. Fahr, C.~Offermanns, A.~Mhamdi, L.~T. Biegler, and A.~Mitsos,
	``Optimal start-up of air separation processes using dynamic optimization
	with complementarity constraints,'' {\em Computer Aided Chemical
		Engineering}, vol.~48, pp.~1147--1152, 2020.
	
	\bibitem{doi:10.7566/JPSJ.92.044802}
	Y.~Furue, M.~Konoshima, H.~Tamura, and J.~Ohkubo, ``Efficient correlation-based
	discretization of continuous variables for annealing machines,'' {\em Journal
		of the Physical Society of Japan}, vol.~92, no.~4, p.~044802, 2023.
	
	\bibitem{Caspari_routing_2022}
	A.~Caspari, S.~Fahr, and A.~Mitsos, ``Optimal eco-routing for hybrid vehicles
	with powertrain model embedded,'' {\em IEEE Transactions on Intelligent
		Transportation Systems}, vol.~23, no.~9, pp.~14632--14648, 2022.
	
\end{thebibliography}

\begin{appendices}
	
	\section{Chemical Reactor Model}
	\label{appendix:chemicalReactorModel}
	
	We use the chemical reactor example from \cite{rawlings2017model}.
	The process model is given by the following system of ordinary differential equations:
	
	\begin{subequations}
		\begin{align}
			\frac{\mathrm{d} c}{\mathrm{d} t} (t) =  &   \frac{F^0(t)(c^0 - c(t))}{\pi r^2 h(t)} \nonumber  \\ &\quad - k^0 \exp\left(-\frac{E}{R T(t)}\right) c(t) 
			\label{eq:pannocchiaReactorNominalModel:c} \\
			\frac{\mathrm{d} T}{\mathrm{d} t} (t) = &  \frac{F^0(t) (T^0 - T(t))}{\pi r^2 h(t)}  \nonumber \\ &\quad -  \frac{\Delta H}{\rho C_\mathrm{p}} k^0 \exp\left(-\frac{E}{R T(t)}\right) c(t) \nonumber \\ &\quad + \frac{2 U}{r \rho C_\mathrm{p}} \left(T^c(t) - T(t)\right) 
			\label{eq:pannocchiaReactorNominalModel:T} \\
			\frac{\mathrm{d} h}{\mathrm{d} t} (t) = &  \frac{F^0(t) - F(t)}{\pi r^2}. \label{eq:pannocchiaReactorNominalModel:h}
		\end{align}
		\label{eq:pannocchiaReactorNominalModel}
	\end{subequations}
	The feed stream enters the reactor at the temperature $T^0$ with the flowrate $F^0$ and the concentration $c^0$. 
	A chemical reaction takes place inside the reactor.
	The reactor is surrounded by a cooling jacket, where cooling water is fed at a temperature of $T^c$.
	$F$ is the outlet flow rate, $c$ is the outlet concentration, $T$ is the outlet temperature, and $h$ is the reactor height.

	To obtain a discrete-time formulation as required for the dynamic optimization problem formulation \eqref{eq:optimizationProblem}, we apply an Euler-forward discretization of the process model \eqref{eq:pannocchiaReactorNominalModel} with equal-spaced time steps of $N = 20$ and $\Delta t = 0.2 $ min.
	This leads to the following discrete-time process model:
	
	\begin{subequations}
		\begin{align*}
			c_{t+1} & =   c_t +  \nonumber  \\
			& \ \ \ \left(\frac{F^0 (c^0 - c_t)}{\pi r^2 h_t}  - k^0 \exp\left(-\frac{E}{R T_t}\right) c_t \right) \cdot \Delta t 
			\\
			T_{t+1} & = T_t + \left( \frac{F^0 (T^0 - T_t)}{\pi r^2 h_t} - \right. \nonumber \\ & \ \ \frac{\Delta H}{\rho C_\mathrm{p}} k^0 \exp\left(-\frac{E}{R T_t}\right) c_t +  \nonumber \\ & \ \  \left. \frac{2 U h_t}{r \rho C_\mathrm{p}} \left(T^c_{t} - T_t\right)  \right) \cdot \Delta t  \\
			h_{t+1} & = h_t + \left( \frac{F^0 - F_t}{\pi r^2} \right) \cdot \Delta t
		\end{align*}
		\label{eq:pannocchiaReactorNominalModel:discrete}
	\end{subequations}

	We assume quasi-stationarity for the variable $h$, leading to $F(t) = F^0(t) = 0.1 \mathrm{m}^3/\mathrm{min}$.
	Furthermore, we use a constant reaction rate and replace $\exp(-\frac{E}{R T(t)})$ by $\exp(-\frac{E}{R T^\mathrm{fix}})$, with $T^\mathrm{fix} = 340 K$, which is required as current quantum annealing hardware can not easily deal with non-polynomial functions.
	These assumptions lead to the discrete-time model \eqref{eq:simplifiedpannochia}.
	In addition, we force the solution to adhere the the bounds $T^c \in [295, 330]$.
	
	Table \ref{tab:caseStudyChemicalReactor:Parameters} lists the values chosen from \cite{rawlings2017model} for the parameters and the value for the reactor height $h$.

	\begin{table}[h]
		\centering
		\caption{Parameter values for chemical reactor case study.}
		\begin{tabular}{ll} \hline
			\textbf{parameter} &\textbf{value} \\ 
			\hline 
			$F^0$& 0.1  m$^3$/min\\ 
			
			$T^0$& 350  K \\ 
			
			$c^0$ & 1 $\mathrm{kmol/m^3}$\\ 
			$k^0$ & 7.2 $\cdot 10^{10}$  min$^{-1}$\\ 
			
			$r$ & 0.219 m\\ 
			$E/R$ & 8750 K\\ 
			$U$ & 54.94 $\mathrm{kJ/min \ m^2 \ K}$\\ 
			$\rho$ & 1000 $\mathrm{kg/m^3}$ \\ 
			$C_\mathrm{p}$ & 0.239 $\mathrm{kJ/ kg \ K}$\\ 
			$\Delta H$ & -5 $\cdot 10^{-4}$ $\mathrm{kJ/kmol}$\\ 
			$h$ & 0.8 $\mathrm{m}$\\ 
			\hline
		\end{tabular} 
		\label{tab:caseStudyChemicalReactor:Parameters}
	\end{table}

	\section{Dynamic optimization with chemical reactor model embedded}
	\label{appendix:DynamicOptimizationProblem}
	
	This section provides the dynamic optimization problem with the chemical reactor model embedded.
	The general form of the the dynamic optimization problem is given by \eqref{eq:optimizationProblem}.
	\autoref{appendix:chemicalReactorModel} presents the chemical reactor model.

	We use the following tracking objective function which penalizes deviations of the reactor outlet temperature $T$ from a target temperature $T^\mathrm{fix}$:
	\begin{equation*}
		\begin{aligned}
			E = \sum_{i=0}^{N}(T_i-T^{\mathrm{}})^2 \, .
		\end{aligned}
	\end{equation*}
	
	We apply the model assumptions as described in \autoref{appendix:chemicalReactorModel}.
	We use a time horizon of 4 min. 
	This is a sufficient time for the CSTR model to reach the target steady state.
	The model is discretized with a time discretization of $\Delta t = 0.2 $ min leading to $N = 20$ equal-spaced time steps.
	
	The resulting discrete-time dynamic optimization problem is given by the following problem statement:
	
	\begin{subequations}
		\begin{align*}
			\min &  \sum_{i=0}^{N=20}(T_i-T^{\mathrm{fix}})^2,  \\
			s.t. \qquad  & c_{i+1}  =  c_{i} + \Delta t \cdot \frac{F^0(c^0 - c_i)}{\pi r^2 h}  \nonumber\\ 
			& \quad - \Delta t \cdot k^0 \exp\left(-\frac{E}{R T^{\mathrm{fix}}}\right) c_{i}, \nonumber \\ & \quad \forall i \in \{0,...,19\}, \\
			& T_{i+1}  = T_i +  \Delta t \cdot \frac{F^0(T^0 - T_i)}{\pi r^2 } \nonumber \\ 
			& \quad -  \Delta t \cdot \frac{\Delta H}{\rho C_\mathrm{p}} k^0 \exp\left(-\frac{E}{R T^{\mathrm{fix}}}\right) c_i \nonumber \\ 
			& \quad  + \Delta t \cdot \frac{2 U}{r \rho C_\mathrm{p}} \left(T^c_i - T_i\right),  \nonumber \\ & \quad \forall i \in \{0,...,19\}, \\
			&  T_0 = 324.5, \\ 
			&  c_0 = 877, \\ 
			& T^c_i \in [295, 330], \forall i \in \{0,...,20\},\\
			&  c_i \in \mathbb{R}, \forall i \in \{1,...,20\}, \\
			&  T_i \in \mathbb{R}, \forall i \in \{1,...,20\}.
		\end{align*}
	\end{subequations}

	\section{Visualization of Problem Embedding on Pegasus Topology}
	\label{appendix:pegasusembedding}
\begin{figure*}[t]
		\centering
		\includegraphics[width=\textwidth]{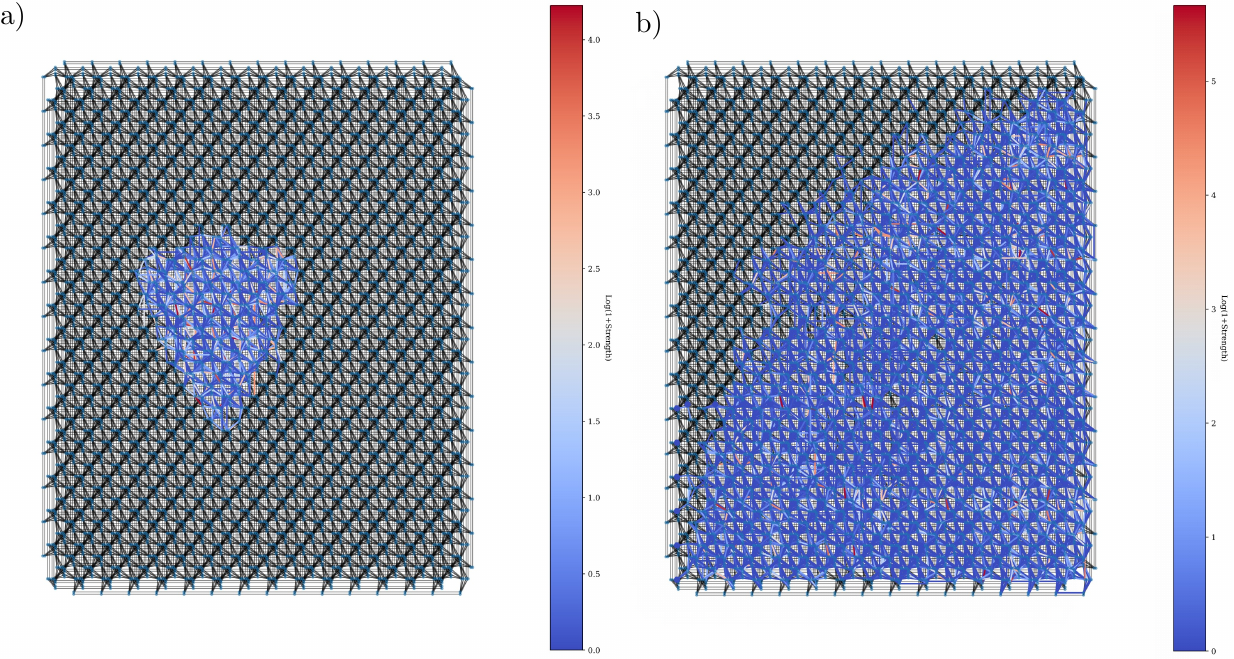}
		\caption{Illustration of the problem graph embedding for two scenarios: (a) a smaller scale (10 timesteps with six binary variables each, totaling approximately 500 interlinked qubits) and (b) a larger scale (20 timesteps with eight binary variables each, cumulating in roughly 3700 interconnected qubits) on the D-Wave Advantage's Pegasus architecture. Here, circles symbolize qubits, while the connecting edges or lines signify the couplers. The coupling strength is indicated by the color of the edges.}    \label{fig:pegasusembedding}
	\end{figure*}

	The D-Wave Advantage QPUs, based on the Pegasus graph topology, represent an evolution from the previous D-Wave 2000Q systems \footnote{D-Wave Systems Inc., "Advantage Processor Overview", \url{https://www.dwavesys.com/media/3xvdipcn/14-1058a-a_advantage_processor_overview.pdf}, accessed: 02-02-2024}. With over 5,000 qubits and increased connectivity—featuring 15 couplers per qubit, i.e. more than 35000 couplers in total—the Advantage QPUs can embed and solve larger and more complex optimization problems than their predecessors. This expansion allows for more intricate mappings of problems like the one expressed in equation~\eqref{eq:qubo_opt}, as visualized in Figure~\ref{fig:pegasusembedding}.
	
\end{appendices}

\end{document}